\documentclass[12pt]{article}
\usepackage{latexsym}
\usepackage{longtable}

\newcommand{\duk}{\noindent {\bf Proof. }}
\newcommand{\kduk}{\hfill $\Box$\bigskip}

\newcommand{\N}{\mathbf{N}}
\newcommand{\Z}{\mathbf{Z}}
\newcommand{\Q}{\mathbf{Q}}
\newcommand{\ex}{\mathrm{ex}}

\newtheorem{veta}{Theorem}[section]
\newtheorem{dusl}[veta]{Corollary}

\newtheorem{prop}[veta]{Proposition}

\def\cla#1#2#3#4#5#6{
  {\sc #1, }#2, {\it #3, }{\bf #4 }(#5), #6.}

\def\pre#1#2#3#4#5{
  {\sc #1, }#2, {\it #3, }technical report {\bf #4}, #5.}

\def\kni#1#2#3#4#5{
  {\sc #1, }{\it #2, }#3, #4, #5.}

\def\vsbo#1#2#3#4#5#6#7#8{
  {\sc #1, }#2. In: {#4 (ed.), } {\it #5, } #6, #7,  #8; pp. #3.}

\begin{document}

\author{Martin Klazar\thanks{Department of Applied Mathematics (KAM) and Institute for Theoretical 
Computer Science (ITI), Charles University, Malostransk\'e n\'am\v est\'\i\ 25, 118 00 Praha, 
Czech Republic. ITI is supported by the project LN00A056 of the 
Ministery of Education of the Czech Republic. E-mail: {\tt klazar@kam.mff.cuni.cz}}}
\title{Extremal problems for ordered hypergraphs: small patterns and some enumeration}
\date{}

\maketitle
\begin{abstract}
We investigate extremal functions $\ex_e(F,n)$ and $\ex_i(F,n)$ counting maximum numbers of edges and maximum 
numbers of vertex-edge incidences in simple hypergraphs $H$ which have $n$ vertices and do not contain a fixed 
hypergraph $F$; the containment respects linear orderings of vertices. We determine both functions exactly if 
$F$ has only distinct singleton edges or if $F$ is one of the 55 hypergraphs with at most four incidences 
(we give proofs only for six cases). We prove some exact formulae and recurrences for the numbers of hypergraphs, 
simple and all, with $n$ incidences and derive rough logarithmic asymptotics of these numbers. Identities analogous to 
Dobi\`nski's formula for Bell numbers are given. 
\end{abstract}

\section{Introduction and definitions}

In this article we consider problems on hypergraphs of the following type. Suppose that $H$ is a simple hypergraph with 
$n$ vertices, which means that $H$ is a finite set of finite nonempty subsets of $\N=\{1,2,\dots\}$ with 
$|\bigcup H|=n$, such that for no three vertices $a<b<c$ in $\bigcup H$ and for no two distinct edges $A$ and $B$ in $H$ 
one has the four incidences $a,b\in A\;\&\;b,c\in B$. What are, in terms of $n$, the maximum possible size $|H|$ and the 
maximum possible number of incidences $\sum_{A\in H}|A|$ of $H$? What are the maxima if the forbidden incidence pattern 
is, for example, 
$a\in A\;\&\;a\in B\;\&\;a,b\in C$ ($a<b$ are vertices and $A,B$, and $C$ are distinct edges)? How many 
distinct hypergraphs with linearly ordered vertices and $n$ incidences, simple and all, are there? The first two questions, 
and quite a few similar ones, are answered in Section 3. The third question is addressed in Section 4. This article is a 
continuation of Klazar \cite{klaz}. We refer the reader to \cite{klaz} for further results and for motivation 
of our extremal problems.   

We denote $\N=\{1,2,3,\dots\}$ and work with the standard linear order $<$ on $\N$. If $a,b,n\in\N$ 
with $a\le b$, we write $[a,b]$ for the interval $\{a,a+1,\dots,b\}$ and $[n]=[1,n]$ for $\{1,2,\dots,n\}$.
A {\em hypergraph\/} $H=(E_i:\ i\in I)$ is a finite list of finite nonempty subsets $E_i$ of 
$\N$, called {\em edges\/}. $H$ is {\em simple\/} if $E_i\neq E_j$ for every $i,j\in I$, 
$i\neq j$. The elements of $\bigcup H=\bigcup_{i\in I} E_i\subset\N$ are called {\em vertices\/}. 
Note that our hypergraphs have no isolated vertices. The {\em simplification\/} of $H$ is the simple hypergraph 
obtained from $H$ by keeping from each family of mutually equal edges just one edge. 
The {\em deletion\/} of $E_j$, $j\in I$, from $H=(E_i:\ i\in I)$ yields the hypergraph $(E_i:\ i\in I')$ where 
$I'=I\backslash\{j\}$. The {\em deletion\/} of $a\in\bigcup H$ from $H$ yields the hypergraph 
$(E_i\backslash\{a\}:\ i\in I)$ where the
$\emptyset$'s arising from $E_i=\{a\}$ are omitted; this operation in general destroys simplicity. We may also delete $a$ 
only from some specified edges. The {\em degree\/} $\deg(v)=\deg_H(v)$ of a vertex $v$ of $H$ is the number of the edges
$E\in H$ such that $v\in E$.
The {\em order\/} $v(H)$ of $H=(E_i:\ i\in I)$ is the number of vertices $v(H)=|\bigcup H|$, the 
{\em size\/} $e(H)$ is the number of edges $e(H)=|H|=|I|$, and the {\em weight\/} $i(H)$ is the number of 
incidences between vertices and edges $i(H)=\sum_{i\in I}|E_i|$. Trivially, $v(H)\le i(H)$ and 
$e(H)\le i(H)$ for every $H$.

Two hypergraphs $H=(E_i:\ i\in I)$ and $H'=(E_i':\ i\in I')$ are {\em isomorphic\/} if there are an 
{\em increasing\/} bijection $F:\bigcup H'\to\bigcup H$ and a bijection $f: I'\rightarrow I$ such that 
$F(E_i')=E_{f(i)}$ for every $i\in I'$. $H'$ is a {\em reduction\/} 
of $H$ if $I'\subset I$ and $E_i'\subset E_i$ for every $i\in I'$. $H'$ is {\em contained\/} in $H$, in symbols 
$H'\prec H$, if $H'$ is isomorphic to a reduction of $H$. We call that reduction of $H$ an $H'$-{\em copy\/} in $H$.
For example, if $H'=(\{1\}_1,\{1\}_2)$ ($H'$ is a singleton edge repeated twice) then $H'\prec H$
if and only if $H$ has two intersecting edges. Another example: If $H'=(\{1,4\},\{2,3\})$ then $H'$ is contained 
in $H$ if and only if $H$ has four vertices $a<b<c<d$ such that $a$ and $d$ lie in one edge of $H$ while $b$ and $c$ 
lie in another edge. If $H'\not\prec H$, we say that $H$ is $H'$-{\em free\/}.
Let $F$ be any hypergraph. We associate with $F$ the extremal functions $\ex_e(F,\cdot),\;\ex_i(F,\cdot): \N\to\N$, 
defined by
\begin{eqnarray*}
\ex_e(F,n)&=&\max\{e(H):\  H\not\succ F\;\&\;\mbox{$H$ is simple}\;\&\;v(H)=n\}\\
\ex_i(F,n)&=&\max\{i(H):\  H\not\succ F\;\&\;\mbox{$H$ is simple}\;\&\;v(H)=n\}.
\end{eqnarray*}
In \cite{klaz} we defined both functions with the requirement $v(H)\le n$. Here we are more interested in their precise 
values and therefore we require $v(H)=n$. 

Obviously, for every $n\in\N$ and 
$F$, $\ex_e(F,n)\le 2^n-1$ and $\ex_i(F,n)\le n2^{n-1}$, but much better bounds can be usually given.
The {\em reversal\/} of a hypergraph $H=(E_i:\ i\in I)$ with $N=\max(\bigcup H)$ is the hypergraph 
$\overline{H}=(\overline{E_i}:\ i\in I)$ where $\overline{E_i}=\{N-x+1:\ x\in E_i\}$. Reversals are 
obtained by reverting the linear ordering of vertices. It is clear that 
$\ex_e(F,n)=\ex_e(\overline{F},n)$ and $\ex_i(F,n)=\ex_i(\overline{F},n)$ for every $F$ and $n$. 

In this article we complement the results of \cite{klaz}, where we derived some asymptotic upper bounds, and 
determine precise values of 
$\ex_e(F,n)$ and $\ex_i(F,n)$ for several hypergraphs $F$. Then we address some naturally arising enumerative questions. 
The present article is a revised version of about one half of the technical report \cite{klaz01}; the other half appears in 
\cite{klaz}. Sections 2 and 3 contain extremal results.
In Theorems~\ref{Ske} and \ref{Ski} we determine $\ex_e(F,n)$ and $\ex_i(F,n)$ exactly if 
$F=S_k=(\{1\},\{2\},\ldots,\{k\})$ consists only of distinct singleton edges. Then both functions are not 
nondecreasing: $\ex_e(S_k,k-1)>\ex_e(S_k,k)$ and $\ex_i(S_k,k-1)>\ex_i(S_k,k)$  ($k\ge 3$). 
In Theorem~\ref{mono} we prove that if $F$ is nonisomorphic to 
$S_k$, then $\ex_e(F,n)<\ex_e(F,n+1)$ for every $n\in\N$. 
Since all hypergraphs obtained from $S_k$ by permuting its vertices are mutually isomorphic, in Theorems~\ref{Ske} and 
\ref{Ski} the ordering of vertices is irrelevant. In Section 3 we determine both extremal functions exactly 
for every of the 55 hypergraphs $F$ with $1\le i(F)\le 4$. In Propositions~\ref{F6}--\ref{F30} we present proofs only for 
six cases (other three cases are subsumed in Theorems~\ref{Ske} and \ref{Ski}). Section 4 is enumerative. 
In Proposition~\ref{hvn} we enumerate simple hypergraphs with order $n$. Theorem~\ref{hranovytyp} enumerates both simple 
and all hypergraphs with prescribed numbers of edges of each cardinality. Corollary~\ref{hnishns} enumerates both simple 
and all hypergraphs with weight $n$ by a sum over integer partitions. Proposition~\ref{efectrecurr} does the same less 
elegantly but more efficiently by 
recurrences. In Corollary~\ref{dobinski} we give identities for hypergraphs which are analogous to the Dobi\`nski's 
formula for set partitions. In Proposition~\ref{asym} we bound the numbers of hypergraphs with weight $n$ by 
the Bell numbers.

\section{Singleton hypergraphs} 

Note that functions $\ex_e((\{1\}),n)$ and $\ex_i((\{1\}),n)$ are undefined.

\begin{veta}\label{Ske}
Let $S_k=(\{1\},\{2\},\ldots,\{k\})$. Then, for $k\ge 2$,
$$
\ex_e(S_k,n)=\left\{
\begin{array}{lll}
2^n-1 & \ldots & 1\le n<k\\
2^{k-2} & \ldots &  n\ge k.
\end{array}
\right.
$$
In particular, for $k\ge 3$ the function $\ex_e(S_k,n)$ has the unique global maximum $\ex_e(S_k,k-1)=2^{k-1}-1$. 
\end{veta}
\duk 
The case $1\le n<k$ is clear. For $n\ge k\ge 2$ we have $\ex_e(S_k,n)\ge 2^{k-2}$ because  
$\{[n]\}\cup(E:\ \emptyset\ne E\subset[k-2])\not\succ S_k$. We prove by induction on $k$ that for $n\ge k$ also  
$\ex_e(S_k,n)\le 2^{k-2}$. For $k=2$ this holds because $\ex_e(S_2,n)=1$ for every $n\in\N$. Let 
$n\ge k\ge 3$ and let $H$ be a simple $S_k$-free hypegraph with $\bigcup H=[n]$. We show that we can assume that (i) 
$\deg(v)\ge 2$ for every $v\in\bigcup H$ and (ii) there is an $E\in H$ with $|E|\ge 2$ and an $a\in E$ such that 
$E\backslash\{a\}\not\in H$. 

If (i) is false, there is a vertex contained in a unique edge. We delete the edge from $H$ and obtain a 
hypergraph $H'$ which must be $S_{k-1}$-free. We are done by induction: 
$e(H)=e(H')+1\le (2^{(k-1)-1}-1)+1=2^{k-2}$. Suppose that (ii) is false. Let $a\in\bigcup H$ be 
arbitrary and $E\in H$, $a\in E$, be such that $|E|$ is as small as possible. If $|E|>1$, there is a $b\in E$, 
$b\neq a$. By the negation of (ii), $E\backslash\{b\}\in H$, contradicting the minimality of $|E|$. 
Thus $|E|=1$ and $\{a\}\in H$. Hence $\{a\}\in H$ for every $a\in\bigcup H$. But this implies the 
contradiction $H\succ S_k$ (since $n\ge k$). 

Thus (i) and (ii) hold. Let $a$ and $E$ be as in (ii). Let $E'\in H$ be such that $a\in E'$, $E'\neq E$, and, 
if possible, $|E'|=1$. We obtain $H'$ by deleting $E'$ from $H$ and then deleting $a$ from $H\backslash\{E'\}$. 
Some edges 
may get duplicated and therefore we set $H''$ to be the simplification of $H'$. By (i), $v(H'')=v(H)-1=n-1\ge k-1$. Since 
any $S_{k-1}$-copy in $H''$ can be extended by $E'$ and $a$ to an $S_k$-copy in $H$, $H''\not\succ S_{k-1}$. 
Also, $e(H')\le 2e(H'')-1$ because, by (ii), $E\backslash\{a\}$ is not duplicated in $H'$. Notice that 
$\emptyset\not\in H''$ because we have deleted $\{a\}$ as $E'$.  By induction (now we use the stronger upper bound 
on $e(H'')$),
$$
e(H)=e(H')+1\le (2e(H'')-1)+1=2e(H'')\le 2\cdot 2^{(k-1)-2}=2^{k-2}. 
$$
\kduk

\noindent
The function $\ex_e(S_k,n)$ has the strange feature of being independent of $n$. We show that other extremal 
functions $\ex_e(F,n)$ are increasing, as one expects. 

\begin{veta}\label{mono}
If $F$ is not isomorphic to any $S_k=(\{1\},\{2\},\dots,\{k\})$, then 
$$
\ex_e(F,n)<\ex_e(F,n+1)
$$ 
for every $n\in\N$. 
\end{veta}
\duk
Let $\bigcup F=[m]$, $m\ge 2$, and $F\ne S_m$. We say that $\{i\}\in F$ is an {\em isolated singleton\/} of $F$ if 
$\deg(i)=1$. Let $l$ be the maximum number such that $\{1\},\{2\},\ldots,\{l\}$ are isolated singletons 
of $F$. Since $F\neq S_m$, we have $0\le l<m$. Any other isolated singleton of $F$ is preceded by at least 
$l+1$ vertices. 
We proceed by induction on $n$. The inequality holds for every $n<m-1$ because then $\ex_e(F,n)=2^n-1$. Let 
$n\ge m-1$ and let $H$, $\bigcup H=[n]$, attain the value $\ex_e(F,n)$. If $a\in E\in H$ and $\{a\}\not\in H$, we replace 
$E$ by $\{a\}$. The new hypergraph is simple, $F$-free, and has the same size as $H$. By the 
inductive assumption, it must have also the same order. Repeating the replacements, 
we obtain a simple $F$-free hypergraph $H'$ such that $e(H')=e(H)=\ex_e(F,n)$, $\bigcup H'=\bigcup H=[n]$, 
and $\{a\}\in H'$ for every $a\in[n]$. We define $H''$ by inserting in $H'$, between the vertices $l$ and $l+1$, a new 
singleton edge $\{u\}$. $H''$ is simple and satisfies $v(H'')=n+1$ and $e(H'')=e(H')+1=\ex_e(F,n)+1$. We show 
that $H''$ is $F$-free. This gives $\ex_e(F,n+1)\ge e(H'')>\ex_e(F,n)$. If $H''\succ F$, the new edge $\{u\}$ 
would have to participate in every $F$-copy in $H''$ as an isolated singleton. It cannot play the role of any of the initial 
$l$ isolated singletons of 
$F$ because $\{i\}\in H'$ for every $i\in[n]$ and $n\ge m-1\ge l$; we would have already $F\prec H'$. It cannot play the role of 
any other isolated singleton of $F$ either because those are preceded in $F$ by at least $l+1$ vertices but 
$\{u\}$ is preceded in $H''$ by only $l$ vertices. Thus $H''\not\succ F$.
\kduk

\begin{veta}\label{Ski}
Let $S_k=(\{1\},\{2\},\ldots,\{k\})$. Then, for $k\ge 2$,
$$
\ex_i(S_k,n)=\left\{
\begin{array}{lll}
n2^{n-1} & \ldots & 1\le n<k\\
n+(k-2)2^{k-3} & \ldots & k\le n\le 2^{k-3}+1\\
(k-1)n-(k-2) & \ldots &  n\ge \max(k,2^{k-3}+1).
\end{array}
\right.
$$
In particular, $\ex_i(S_k,k-1)>\ex_i(S_k,n)$ for $k\le n\le \max(k,2^{k-2})$ ($k\ge 3$). 
\end{veta}
\duk
The first case is clear. We suppose that $n\ge k\ge 2$ and that $H$ is a simple hypergraph with 
$\bigcup H=[n]$. We consider its dual $H^*$:
$$
H^*=(E^*_i:\ i\in[n])\mbox{ where }E^*_i=\{E\in H:\ i\in E\}.
$$
Thus $e(H^*)=v(H)=n$. Let $\Gamma(X)=\Gamma_{H}(X)$ be for 
$X\subset[n]$ defined by
$$
\Gamma(X)=\left|\bigcup_{i\in X}E^*_i\right|=
|\{E\in H:\ E\cap X\neq\emptyset\}|.
$$
By the defect form of P. Hall's theorem (Lov\'asz \cite[Problems 7.5 and 13.5]{lova}) applied on $H^*$, $H$ 
is $S_k$-free if and only if 
$$
\max_{X\subset[n]}|X|-\Gamma(X)\ge n-k+1.
$$
Thus if $H$ is $S_k$-free, there exists a set $X\subset[n]$ of cardinality $l$, $n-k+2\le l\le n$ 
($\Gamma(X)\ge 1$), intersected by only at most $l-n+k-1$ edges of $H$. And contrarywise, every such a hypergraph 
is (trivially) $S_k$-free. Hence
$$
i(H)\le (l-n+k-1)n-(l-n+k-2)+(n-l)2^{n-l-1}=f(l,k,n)
$$
and this bound is attained.

Let $k$ and $n$ be fixed. The first difference of $f(l,k,n)$ with respect to $l$ is the increasing function 
$$
f(l+1,k,n)-f(l,k,n)=n-1-(n-l+1)2^{n-l-2}.
$$ 
Therefore $f(l,k,n)$ attains its maximum in one of the endpoints $l=n-k+2$ and $l=n$ or in both. The 
corresponding values are $f(n-k+2,k,n)=n+(k-2)2^{k-3}$ and $f(n,k,n)=(k-1)n-(k-2)$. These values 
are equal for $n=2^{k-3}+1$. For $n<2^{k-3}+1$ the former value dominates and for $n>2^{k-3}+1$ the latter. We 
obtain the values of $\ex_i(S_k,n)$ in the remaining two cases. Maximum weights are attained by the hypergraph $H_1$ or 
by $H_2$, where the edges of $H_1$, respectively of $H_2$, are $[n]$ and all nonempty subsets 
of some $(k-2)$-element set $Y\subset[n]$, respectively $[n]$ and some $k-2$ distinct $(n-1)$-element 
subsets of $[n]$.  
\kduk

\noindent
For $1\le n<k$ the maximum weight is attained only by the complete hypergraph.
The proof shows that for $n\ge k$ the only types of extremal hypergraphs are $H_1$ and $H_2$. Thus the number of 
simple $S_k$-free hypergraphs having order $n$ and the maximum weight equals $1$ 
if $1\le n<k$ and equals $\eta{n\choose k-2}$ if $n\ge k$, where for $k=2,3,4$ 
always $\eta=1$ and for $k\ge 5$ we have $\eta=1$ if 
$n\neq 2^{k-3}+1$ and $\eta=2$ if $n=2^{k-3}+1$.

One can use P. Hall's theorem to give another proof of Theorem~\ref{Ske}. The number of hypergraphs $H$ attaining 
the value $\ex_e(S_k,n)$ is seen to be $1$ for $n<k$ and $2^{k-2}{n\choose k-2}$ for $n\ge k$. The latter hypergraphs are all 
$H$ of the form $H=\{Y\}\cup(X:\ \emptyset\ne X\subset Z)$ where $Z$ is a $k-2$-element subset of $[n]$ and 
$[n]\backslash Z\subset Y\subset[n]$.

We conjecture that if $F$ is not isomorphic to any of the singleton hypergraphs $S_k=(\{1\},\{2\},\dots,\{k\})$, then
$$
\ex_i(F,n)<\ex_i(F,n+1)
$$
for every $n\in\N$. 

\section{Forbidden hypergraphs of weight at most $4$}

In this section we give precise formulae for $\ex_e(F,n)$ and $\ex_i(F,n)$ for every $F$ with $1\le i(F)\le 4$. 
There are 
55 such nonisomorphic hypergraphs but due to the reversals it suffices to consider 39 of them. The proofs are 
usually straightforward and often repetitive. Lest the reader be bored and tired, we present here only a sample consisting 
of six cases. The proofs for all of the 39 cases can be found in \cite{klaz01}. First we list the hypergraphs $F$, then 
we review the results in a table, and in the rest of the section we give proofs for six cases.

Weight 1 and 2: 
$$
F_1=(\{1\}),\ \ F_2=(\{1\}_1,\{1\}_2),\ \ F_3=(\{1\},\{2\}),\ \mbox{ and }\ F_4=(\{1,2\}).
$$
Weight 3: 
$$
F_5=(\{1\}_1,\{1\}_2,\{1\}_3),\ \ F_6=(\{1\}_1,\{1\}_2,\{2\}),\ \ \overline{F_6},\ \ F_7=(\{1\},\{2\},\{3\}), 
$$
$$
F_8=(\{1\},\{1,2\}),\ \ \overline{F_8},\ \ F_9=(\{1\},\{2,3\}),\ \ \overline{F_9},\ \ F_{10}=(\{1,3\},\{2\}),
$$
and 
$$
F_{11}=(\{1,2,3\}).
$$
Weight 4:
$$
F_{12}=(\{1\}_1,\{1\}_2,\{1\}_3,\{1\}_4),\ \ F_{13}=(\{1\}_1,\{1\}_2,\{1\}_3,\{2\}),\ \ \overline{F_{13}},
$$
$$
F_{14}=(\{1\}_1,\{1\}_2,\{2\}_1,\{2\}_2),\ \ F_{15}=(\{1\}_1,\{1\}_2,\{2\},\{3\}),\ \ \overline{F_{15}},
$$
$$
F_{16}=(\{1\},\{2\}_1,\{2\}_2,\{3\}),\ \ F_{17}=(\{1\},\{2\},\{3\},\{4\}),
$$
$$
F_{18}=(\{1\}_1,\{1\}_2,\{1,2\}),\ \ \overline{F_{18}},\ \ F_{19}=(\{1\}_1,\{1\}_2,\{2,3\}),\ \ \overline{F_{19}},
$$
$$
F_{20}=(\{1,3\},\{2\}_1,\{2\}_2),\ \ F_{21}=(\{1\},\{2\},\{2,3\}),\ \ \overline{F_{21}},
$$
$$
F_{22}=(\{1\},\{2,3\},\{3\}),\ \ \overline{F_{22}},\ \ F_{23}=(\{1\},\{2\},\{1,3\}),\ \ \overline{F_{23}},
$$
$$
F_{24}=(\{1\},\{2\},\{1,2\}),\ \ F_{25}=(\{1\},\{2\},\{3,4\}),\ \ \overline{F_{25}},
$$
$$
F_{26}=(\{1\},\{2,4\},\{3\}),\ \ \overline{F_{26}},\ \ F_{27}=(\{1\},\{2,3\},\{4\}),
$$
$$
F_{28}=(\{1,4\},\{2\},\{3\}),\ \ F_{29}=(\{1,2\},\{1,3\}),\ \ \overline{F_{29}},
$$
$$
F_{30}=(\{1,2\},\{2,3\}),\ \ F_{31}=(\{1,2\}_1,\{1,2\}_2),\ \ F_{32}=(\{1,2\},\{3,4\}),
$$
$$
F_{33}=(\{1,4\},\{2,3\}),\ \ F_{34}=(\{1,3\},\{2,4\}),\ \ F_{35}=(\{1\},\{1,2,3\}),
$$
$$
\overline{F_{35}},\ \ F_{36}=(\{1,2,3\},\{2\}),\ \ F_{37}=(\{1\},\{2,3,4\}),\ \ \overline{F_{37}},
$$
$$
F_{38}=(\{1,3,4\},\{2\}),\ \ \overline{F_{38}},\ \mbox{ and }\  F_{39}=(\{1,2,3,4\}).
$$

The formulae in the table below hold for every $n\in\N$ if it is not written else. The omitted values are:
$\ex_e(F_k,1)=\ex_i(F_k,1)=1$ for every $k$, $\ex_i(F_7,2)=4$, $\ex_e(F_{12},2)=3$,
$\ex_i(F_{12},2)=4$, $\ex_i(F_{17},3)=12$, $\ex_i(F_{18},2)=4$, $\ex_i(F_{18},3)=8$, $\ex_i(F_{18},4)=11$, 
$\ex_i(F_{18},5)=15$, and $\ex_i(F_{30},3)=8$. In the first column, numbers $\overline{k}$ with bar indicate 
that $F_k$ is nonisomorphic to $\overline{F_k}$ and thus the formulae in the $k$-th row apply to two hypergraphs.

\setlongtables
\begin{longtable}{|l|c|c|}
\hline
\mbox{$k$}& $\ex_e(F_k,n)$ & $\ex_i(F_k,n)$\\
\hline
\endhead
1 &  not defined & not defined\\
2 &  $n$ & $n$\\
3 &  $1,1,\ldots$ & $n$\\
4 & $n$ & $n$ \\
5 & $\left\lfloor 3n/2\right\rfloor$ & $2n\ (n>1)$\\
$\overline{6}$ & $n$ & $2n-1$\\
7 & $1,3,2,2,\ldots$ & $2n-1\ (n\neq 2)$\\
$\overline{8}$ & $n$ & $2n-1$\\
$\overline{9}$ & $2n-1$ & $3n-2$\\
10 & $2n-1$ & $3n-2$ \\
11 & $(n^2+n)/2$ & $n^2$\\
12 & $2n\ (n>2)$ & $3n\ (n>2)$\\
$\overline{13}$ & $2n-1$ & $\left\lfloor 7(n-1)/2\right\rfloor+1$\\
14 & $n+1\ (n>1)$ & $3n-2$\\
$\overline{15}$ & $n+1\ (n>1)$ & $3n-2$\\
16 & $n+1\ (n>1)$ & $3n-2$\\
17 & $1,3,7,4,4,\ldots$ & $3n-2\ (n\neq 3)$\\
$\overline{18}$ & $2n-1$ & $4n-6\ (n>5)$ \\
$\overline{19}$ & $2n-1$ & $3n-2$\\
20 & $2n-1$ & $3n-2$\\
$\overline{21}$ & $2n-1$ & $3n-2$\\
$\overline{22}$ & $2n-1$ & $3n-2$\\
$\overline{23}$ & $2n-1$ & $3n-2$\\
24 & $n$ & $2n-1$\\
$\overline{25}$ & $4n-5\ (n>1)$ & $8n-12\ (n>1)$\\
$\overline{26}$ & $4n-5\ (n>1)$ & $8n-12\ (n>1)$\\
27 & $4n-5\ (n>1)$ & $8n-12\ (n>1)$\\
28 & $4n-5\ (n>1)$ & $8n-12\ (n>1)$\\
$\overline{29}$ & $2n-1$ & $4n-4\ (n>1)$\\
30 & $\left\lfloor n^2/4\right\rfloor+n$ & $2\left\lfloor n^2/4\right\rfloor+n\ (n\neq 3)$\\
31 & $(n^2+n)/2$ & $n^2$\\
32 & $2\left\lfloor(n+1)^2/4\right\rfloor-1$ & $5\left\lfloor(n+1)^2/4\right\rfloor-2n-2$\\
33 & $4n-5\ (n>1)$ & $8n-12\ (n>1)$\\
34 & $4n-5\ (n>1)$ & $8n-12\ (n>1)$\\
$\overline{35}$ & $(n^2+n)/2$ & $n^2$\\
36 & $(n^2+n)/2$ & $n^2$\\
$\overline{37}$ & $n^2-n+1$ & $(5n^2-9n+6)/2$\\
$\overline{38}$ & $n^2-n+1$ & $(5n^2-9n+6)/2$\\
39 & $(n^3+5n)/6$ & $(n^3-n^2+2n)/2$\\
\hline
\end{longtable}

\noindent
The results for $k=3,7$, and $17$ are particular cases of Theorems~\ref{Ske} and \ref{Ski}. Cases $k=33$ and $34$ 
were proved already in Klazar \cite{klaz00}. 

Suppose $H$ is a simple 
hypergraph such that $H\not\succ F$ for some $F$, $E\in H$ is an edge, and $a\in E$ is a vetex such that 
$\{a\}\not\in H$. Replacing $E$ with $\{a\}$ we obtain a hypergraph $H'$ with the same size as $H$ and possibly smaller 
order. Moreover, $H'$ is simple and $H'\not\succ F$. Repeating the replacements, in the end we obtain a {\em 
singleton completion\/} $H'$ of $H$ with these properties: $H'$ is simple, $H'\not\succ F$, $e(H')=e(H)$, 
$v(H')\le v(H)$, and $\{a\}\in H'$ for every $a\in\bigcup H'$. Singleton completion helps to determine $\ex_e(F,n)$ if 
$F$ has at least one singleton edge; we used it already in the proof of Theorem~\ref{mono}. 

\begin{prop}\label{F6}
For every $n\in\N$, $\ex_e(F_6,n)=n$ and $\ex_i(F_6,n)=2n-1$. 
\end{prop}
\duk
We have $\ex_e(F_6,n)\ge n$ because $(\{i\}:\ i\in[n])\not\succ F_6$. Let $H$ be any simple hypergraph 
with $H\not\succ F_6$ and $v(H)=n$ and let $H'$ be its singleton completion, $v(H')=m\le n$. $H'$ 
has no nonsingleton edges and thus $e(H)=e(H')=m\le n$. 

We have $\ex_i(F_6,n)\ge 2n-1$ because $([n],[n-1])\not\succ F_6$. Also, because 
$(\{i,n\},\{n\}:\ i\in [n-1])\not\succ F_6$. Let $H$ be any simple hypergraph with $H\not\succ F_6$ 
and $\bigcup H=[n]$. Then $\deg(a)\le 2$ for every $a\in [n-1]$, and the equality for 
some $a$ implies $\deg(n)\le 2$. Hence $\deg(a)=2$ for an $a<n$ implies $i(H)\le 2n$. In fact, even $i(H)\le 2n-1$ 
because $\deg(a)=2$ for every $a\in[n]$ is impossible ($H$ is simple). In the other case when $\deg(a)=1$ for every 
$a<n$ again $i(H)\le 2n-1$ because then $\deg(n)\le n$. In both cases $i(H)\le 2n-1$.
\kduk

\begin{prop}\label{F5}
For every $n\in\N$, $\ex_e(F_5,n)=\left\lfloor 3n/2\right\rfloor$ and $\ex_i(F_5,n)=2n$ ($n>1$). For every $n>2$, 
$\ex_e(F_{12},n)=2n$ and $\ex_i(F_{12},n)=3n$. 
\end{prop}
\duk
The conditions $H\not\succ F_5$ and $H\not\succ F_{12}$ are equivalent, respectively, with $\deg_H(v)\le 2$ and 
$\deg_H(v)\le 3$ for every $v\in\bigcup H$. Thus the results for $\ex_i(F_5,n)$ and $\ex_i(F_{12},n)$ are clear.

We have $\ex_e(F_5,n)\ge n+\lfloor n/2\rfloor$ because 
$(\{i\},\{2j-1,2j\}:\ i\in[n],j\in[\lfloor n/2\rfloor])\not\succ F_5$. Let $H$ be any simple hypergraph 
with $H\not\succ F_5$ and $v(H)=n$ and let $H'$ be its singleton completion, $v(H')=m\le n$. It follows that  
$e(H)=e(H')\le m+\lfloor m/2\rfloor\le n+\lfloor n/2\rfloor$ because the nonsingleton edges of $H'$ must be mutually 
disjoint. 

We have $\ex_e(F_{12},n)\ge 2n$ ($n>2$) because $(\{i\},\{i,i+1\}\ (\mathrm{mod}\ n):\ i\in[n])\not\succ F_{12}$. 
Let $H$ be any simple hypergraph with $H\not\succ F_{12}$ and $v(H)=n$ and let $H'$ be its singleton completion. 
If $|E|\ge 3$ for an edge $E\in H'$, then $E_1\not\in H'$ for some $E_1\subset E$ with $|E_1|=2$. 
Replacing, one by one, $E$ with $E_1$, we get rid of all edges with three and more vertices. We obtain a simple 
$H''$ such that $H''\not\succ F_{12}$, $v(H'')=m\le n$, $e(H'')=e(H')=e(H)$, $|E|\le 2$ for every $E\in H''$, and 
$\{a\}\in H''$ for every $a\in\bigcup H'$. Hence $e(H)=e(H'')\le m+m\le 2n$ because the 2-element edges of $H''$ 
must form disjoint paths and cycles (every vertex is contained in at most two 2-element edges). 
\kduk

The next result answers our second initial question.
\begin{prop}
For every $n\in\N$, $\ex_e(F_{18},n)=2n-1$. As for the other function, $\ex_i(F_{18},1)=1$, $\ex_i(F_{18},2)=4$, 
$\ex_i(F_{18},3)=8$, $\ex_i(F_{18},4)=11$, $\ex_i(F_{18},5)=15$, and $\ex_i(F_{18},n)=4n-6$ for $n\ge 6$. 
\end{prop}
\duk
We have $\ex_e(F_{18},n)\ge 2n-1$ because $(\{i\},\{i,n\},\{n\}:\ i\in[n-1])\not\succ F_{18}$. 
Let $H$ be any simple hypergraph with $H\not\succ F_{18}$ and $v(H)=n$ and let $H'$ be its singleton completion, 
$v(H')=m\le n$. In $H'$, every two nonsingleton edges may intersect only in the common last vertex. Deleting from each 
nonsingleton edge of $H'$ its last vertex, we obtain mutually disjoint subsets of $[n-1]$. Hence 
$e(H)=e(H')\le m+n-1\le 2n-1$. 

We determine $\ex_i(F_{18},n)$; this is not as easy as it might seem. We have $\ex_i(F_{18},n)\ge 4n-6$ for $n\ge 6$
because $(\{i,n-1\},\{i,n\},\{n-1\},\{n\}:\ i\in[n-2])\not\succ F_{18}$. To prove the opposite inequality, 
consider a simple $F_{18}$-free $H$ with 
$\bigcup H=[n]$. Since $H\not\succ F_{18}$, $\deg(1)\le 2$. We delete $1$ from $H$ and obtain $H_1$; $i(H_1)\le i(H)+2$. 
$H_1$ has at most two duplicated edges. Let $E_1=E_2$ be one of the duplications. If $|E_1|=1$, we delete $E_1$ from $H_1$. If 
$|E_1|\ge 2$, we delete from $E_1$ its last vertex. This creates no new duplication (else $H\succ F_{18}$). In this way 
we remove from $H_1$ both possible duplications and obtain a simple $H_2$ with $\bigcup H_2=[2,n]$ and $i(H)\le 4+i(H_2)$.
We have the inductive inequality $i(H)\le 4+\ex_i(F_{18},n-1)$. Note that $\deg_H(2)\le 2$ and thus for induction we may as 
well delete $2$ instead of $1$. If one of $\{1\}$, $\{2\}$, and $\{1,2\}$ is an edge of $H$, then the deletion of 
$\{1\}$ or $\{2\}$ and the removal of at most one duplication give us the 
stronger bound $i(H)\le 3+\ex_i(F_{18},n-1)$. Note also that $\deg_H(v)\ge 3$ implies that $v$ is the last vertex of every 
edge containing it.

We prove that for $n=1,2,3,4,5,$ and $6$ one has $\ex_i(F_{18},n)=1,4,8,11,15,$ and $18$, and that $\ex_i(F_{18},n)\le 4n-6$ for 
$n\ge 6$. 
The first two values are trivial. By the inductive inequality, $\ex_i(F_{18},3)\le 4+4=8$. Weight $8$ is attained 
by $(\{3\},\{1,3\},\{2,3\},[3])$. Let $n=4$, $H$ be simple and $F_{18}$-free, and $\bigcup H=[4]$. Clearly, 
$\deg(1),\deg(2)\le 2$. Let first 
$\deg(3)\ge 3$ and $p$ be the number of edges in $H$ intersecting both $[2]$ and $[3,4]$. Clearly, $p\le\deg(1)+\deg(2)\le 4$. 
Since no edge can contain both $3$ and $4$, $\deg(3)+\deg(4)\le p+2\le 6$ and $i(H)=\sum_1^4\deg(i)\le 2\cdot 2+6=10$. Now let 
$\deg(3)\le 2$ and $p$ be the number of edges $E\in H$ such that $4\in E$ and $E\cap[3]\neq\emptyset$. Then
$p\le \ex_e(F_5,3)=4$, $\deg(4)\le 1+p\le 5$, and $i(H)=\sum_1^4\deg(i)\le 3\cdot 2+5=11$. Weight $11$ is attained by 
$(\{4\},\{i,4\},[4]:\ i\in[3])$. Thus $\ex_i(F_{18},4)=11$. By the inductive inequality, $\ex_i(F_{18},5)\le 4+11=15$.
Weight $15$ is attained by $(\{5\},\{i,5\},\{2j-1,2j,5\}:\ i\in[4],j\in[2])$.

It remains to show that $\ex_i(F_{18},6)=18$ and not $4+15=19$. Weight 18 is attained by 
$(\{6\},\{i,6\},\{1,2,6\},\{3,4,5,6\}:\ i\in[5])$. We elaborate the argument that we used for $n=4$. 
Let $H$, $\bigcup H=[6]$, be simple and $F_{18}$-free. Clearly, 
$\deg(1),\deg(2)\le 2$ and $\deg(3)\le 4$. If $\deg(3)=4$, no edge intersects both $[3]$ and $[4,6]$ 
and $i(H)\le 2\cdot\ex_i(F_{18},3)=16$. If $\deg(3)=3$, we delete $3$ from $H$. If this creates a duplication, one of 
$\{1\}$, $\{2\}$ or $\{1,2\}$ is an edge of $H$ and by the above remark, $i(H)\le 3+\ex_i(F_{18},5)=18$. If no 
duplication arises, again $i(H)\le\deg(3)+\ex_i(F_{18},5)=18$. So $\deg(3)\le 2$. Let $k=\deg(4)$. Let first $k\ge 3$ 
and $p$ be the number of edges intersecting both $[4]$ and $[5,6]$ (none of them contains $4$). If $4\in E\in H$ then $4=\max E$.
Therefore edges incident with $4$ contribute by at least $k-1$ to $\deg(1)+\deg(2)+\deg(3)\le 6$ and thus $k\le 7$ and 
$p\le 6-(k-1)=7-k$. If $\deg(5)\ge 3$, $\deg(5)+\deg(6)\le p+2\le 9-k$ (no edge contains both $5$ and $6$) and 
$i(H)=\sum_1^6\deg(i)\le 3\cdot 2+k+9-k=15$. If $\deg(5)\le 2$, we have $\deg(6)\le 2+p\le 9-k$ and 
$i(H)\le 4\cdot 2+k+9-k=17$. We may assume that $k=\deg(4)\le 2$ and thus $\deg(i)\le 2$ for every $i\in[4]$. If $\deg(5)\ge 3$, we 
again set $p$ to be the number of edges $E\in H$ intersecting both $[4]$ and $[5,6]$. We have $p\le 4\cdot 2=8$ and 
$\deg(5)+\deg(6)\le p+2\le 10$. Thus $i(H)=\sum_1^6\deg(i)\le 4\cdot 2+10=18$. If $\deg(5)\le 2$, let $p$ be 
the number of edges $E\in H$ intersecting $[5]$ and containing $6$. Then $p\le \ex_e(F_5,5)=7$ and 
$\deg(6)\le 1+p\le 8$. We have again $i(H)=\sum_1^6\deg(i)\le 5\cdot 2+8=18$. Thus $\ex_i(F_{18},6)=18$. 

Finally, using induction starting at $n=6$ and the inductive inequality, we see that for $n\ge 6$ we have 
$\ex_i(F_{18},n)\le 4n-6$.  
\kduk

\noindent
The irregular initial behaviour of $\ex_i(F_{18},n)$ permits to start the induction only from $n=6$. This makes 
$\ex_i(F_{18},n)$ the hardest function of the table to determine.   

We have chosen to present the following case because its treatment in \cite{klaz01} contains errors.
\begin{prop}
For every $n\in\N$, $\ex_e(F_{29},n)=2n-1$. For every $n>1$, $\ex_i(F_{29},n)=4n-4$ (and $\ex_i(F_{29},1)=1$).
\end{prop}
\duk
We have $\ex_e(F_{29},n)\ge 2n-1$ because $(\{i\},\{n\},\{i,n\}:\ i\in[n-1])\not\succ F_{29}$.
Let $H$ be any simple $F_{29}$-free hypergraph with $v(H)=n$. It follows that the first vertices of the nonsingleton 
edges of $H$ must be all distinct. Thus $e(H)\le n+n-1=2n-1$. 

We have $\ex_i(F_{29},n)\ge 4n-4$ (for $n>1$) because $(\{i\},\{j,n-1,n\},\{n-1,n\}:\ i\in[n],j\in[n-2])\not\succ F_{29}$.
Let $H$ be any simple $F_{29}$-free hypergraph with $\bigcup H=[n]$, $n>1$. We delete $1$ from $H$ and obtain $H'$. 
From the previous argument we know that $\deg_H(1)\le 2$. Thus $i(H)\le i(H')+2$. One duplication may appear in $H'$ 
if $A\in H$ and $\{1\}\cup A\in H$ for some $A\subset[2,n]$. If this happens, we delete (one) $A$ from $H'$ and 
obtain $H''$. Else we set $H''=H'$. $H''$ is simple, $F_{29}$-free and $\bigcup H''=[2,n]$. $|A|\ge 3$ implies $H\succ F_{29}$ 
which is forbidden. Thus $|A|\le 2$ 
and we have the inductive inequality $i(H)\le i(H')+2\le i(H'')+4\le\ex_i(F_{29},n-1)+4$. Starting from 
$\ex_i(F_{29},2)=4$, induction shows that $\ex_i(F_{29},n)\le 4n-4$.
\kduk

The next result answers our first initial question.
\begin{prop}\label{F30}
For every $n\in\N$, $\ex_e(F_{30},n)=\left\lfloor n^2/4\right\rfloor+n$. We have  
$\ex_i(F_{30},n)=2\left\lfloor n^2/4\right\rfloor+n$ for $n\ne 3$ and $\ex_i(F_{30},3)=8$.
\end{prop}
\duk
We have $\ex_e(F_{30},n)\ge \left\lfloor n^2/4\right\rfloor+n$ because 
$B_n=(\{i\},\{j,k\}: i\in[n], j\in[\lfloor n/2\rfloor], k\in[\lfloor n/2\rfloor+1,n])\not\succ F_{30}$. 
Let $H$ be any simple $F_{30}$-free hypergraph with $\bigcup H=[n]$. If $|E|\ge 3$ for some $E\in H$, we replace $E$ with 
the two-element set consisting of the first two vertices of $E$. The resulting hypergraph is $F_{30}$-free and,
since $H\not\succ F_{30}$, it is simple. Repeating the replacements, we get rid of all edges with three and more elements and may 
assume that $|E|\le 2$ for every $E\in H$. The two-element edges of $H$ form a triangle-free graph on at 
most $n$ vertices. By a special case of Tur\'an's theorem (see \cite[Problem 10.30]{lova}), 
$e(H)\le n+\lfloor{n^2\over 4}\rfloor$.

The lower bound on $\ex_i(F_{30},n)$ is provided again by $B_n$. We show that the maximum weight is attained 
also by $B_n$ with the exception of $n=3$ when $\ex_i(F_{30},3)=8$ and not $7$. We take any simple $F_{30}$-free 
hypergraph $H$ with $\bigcup H=[n]$ and eliminate large edges. 
If $E=\{a_1,a_2,\ldots,a_t\}\in H$ with $t\ge 4$ and $a_1<a_2<\ldots<a_t$, we replace $E$ 
with the edges $\{a_1,a_{t-1}\}, \{a_2,a_{t-1}\},\ldots, \{a_{t-2},a_{t-1}\}$. The resulting hypergraph $H'$ is simple, 
$F_{30}$-free, and satisfies $v(H')\le v(H)$ and $i(H')\ge i(H)$. In this way we eliminate all edges with four or more 
elements. If $t=3$ and $a_3<n$, we replace $E$ with $\{a_2,a_3\}$ and $\{a_2,n\}$. Similarly if $1<a_1$. Thus for 
bounding $i(H)$ from above we may assume that $|E|\le 3$ for every $E\in H$ and that every 3-element edge, say $H$ 
has $k$ of them, is of the form $\{1,a,n\}$. No two-element edge is incident with any of the $a$'s and they form a 
triangle-free graph on at most $n-k$ vertices. By Tur\'an's theorem, $i(H)\le n+2\lfloor{(n-k)^2\over 4}\rfloor+3k$ 
and the bound is attained. For $n\ge 4$ it is maximized for $k=0$ and for $n=3$ for $k=1$. Indeed, 
$(\{1\},\{2\},\{3\},\{1,3\},\{1,2,3\})$ has weight 8 and $(\{1\},\{2\},\{3\},\{1,2\},\{1,3\})$ has weight 7. 
\kduk

For each $F$ with $i(F)\le 4$ it was not too hard to determine its extremal functions but for $i(F)=5$ or 
$6$ difficult cases start to appear. For example, it would be interesting to know what are $\ex_e(F,n)$ and $\ex_i(F,n)$, or
even the graph version of $\ex_e(F,n)$, if 
$F=(\{1,6\},\{2,5\},\{3,4\})$ or if $F=(\{1,2\},\{2,4\},\{3,5\})$ or if $F$ is some other ordered graph with three edges 
(there are 75 of them, 62 simple, see the table in the next section).   

\section{Enumeration of hypergraphs}

For a hypergraph $F$ and $n\in\N$, we let $h_n(F)$ denote the number of all simple nonisomorphic $F$-free hypergraphs 
$H$ with $v(H)=n$. Let $h_n'(F)$ and 
$h_n''(F)$ be the analogous counting functions with $v(H)=n$ replaced by $i(H)=n$ and with the simplicity 
of $H$ dropped in $h_n''(F)$. Remember that we work with the ordered isomorphism; e.g.,  
$F_{29}=(\{1,2\},\{1,3\})$ and $F_{30}=(\{1,2\},\{2,3\})$ are nonisomorphic. The enumerative problems to determine 
or to bound these counting functions are already for $i(F)\le 4$ much more difficult than the 
extremal problems. It suffices to note, for example, that if $F=F_2=(\{1\}_1,\{1\}_2)$ then 
$h_n(F)=h_n'(F)=h_n''(F)=b_n$ where $b_n$ is the Bell number that counts the partitions of $[n]$. 

In Klazar \cite{klaz00} we found the ordinary generating functions $G_1(x)$, $G_2(x)$, and 
$G_3(x)$ of $h_n(F_{34})$, 
$h_n'(F_{34})$, and $h_n''(F_{34})$, respectively. (Recall that $F_{34}=(\{1,3\},\{2,4\})$.) $G_1$, $G_2$, 
and $G_3$ are algebraic over $\Z(x)$ of degrees 3, 4, and 4, respectively, and their coefficients grow roughly like 
$(63.97055\ldots)^n$, $(5.79950\ldots)^n$, and $(6.06688\ldots)^n$ where the bases of the 
exponentials are algebraic numbers of degrees 4, 15, and 23, respectively. We did not succeed in enumerating 
$F_{33}$-free hypergraphs ($F_{33}=(\{1,4\},\{2,3\})$) and we think it is a problem 
that deserves interest. 

Here we shall investigate the total numbers $h_n$, $h_n'$, and $h_n''$ of, respectively, 
all simple nonisomorphic hypergraphs with $n$ vertices, all simple nonisomorphic hypergraphs with weight $n$, and 
all nonisomorphic hypergraphs with weight $n$. The numbers $h_n$ have been 
considered before in the problem of set covers but the remaining two problems seem 
new. We review the known formulae for $h_n$, derive for them a new recurrence, and then proceed to 
$h_n'$ and $h_n''$. 

\begin{prop}\label{hvn}
The numbers $h_n$ of nonisomorphic simple hypergraphs with $n$ vertices satisfy for every $n\ge 1$ the following 
formulae.
\begin{eqnarray*}
1. \hspace{1cm}h_n&=&2^{2^n-1}-\sum_{j=0}^{n-1}{n\choose j}h_j\label{recprohn}\ \ (h_0=1)\\
2. \hspace{1cm}h_n&=&\sum_{j=0}^n(-1)^{n-j}{n\choose j}2^{2^j-1}\label{explrovprohn}\\
3. \hspace{1cm}h_n&=&2\sum_{k,l\ge 0}{h_kh_l\cdot(n-1)!\over (k+l-n+1)!\cdot(n-1-k)!\cdot(n-1-l)!}-h_{n-1}
\end{eqnarray*}
where in 3 the summation range is $\max(k,l)\le n-1\le k+l$.
\end{prop}
\duk
1. This recurrence is proved in Hearne and Wagner \cite{hear_wagn} and is a rearrangement of the identity 
$$
2^{2^n-1}=\sum_{j=0}^n{n\choose j}h_j.
$$
The identity follows by noting that every simple hypergraph with $j\le n$ vertices is isomorphic to exactly 
${n\choose j}$ hypergraphs $H$ with $v(H)=j$ and $\bigcup H\subset[n]$, and that the simple hypergraphs $H$ with 
$\bigcup H\subset[n]$ correspond bijectively to the elements of the power set of the set $\{X\subset[n]:\ X\ne\emptyset\}$.

2. This formula is proved in Comtet \cite[p. 165]{comt} and also in Macula \cite{macu}. We note that the identity of 1 
is equivalent to $F(x)=\mathrm{e}^xH(x)$ where 
$$
F(x)=\sum_{n\ge 0}{2^{2^n-1}x^n\over n!}\ \mbox{ and }\ 
H(x)=\sum_{n\ge 0}{h_nx^n\over n!}
$$
are exponential generating functions of the involved quantitites. Thus $H(x)=\mathrm{e}^{-x}F(x)$ and 
the formula follows.

3. This recurrence follows from the combinatorial definition of $h_n$. Any simple hypergraph $H$ with 
$\bigcup H=[n]$ decomposes uniquely into two hypergraphs $H_1$ and $H_2$: $H_1$ consists of the sets $E\backslash\{1\}$ 
such that $1\in E\in H$ (we omit the $\emptyset$ if $\{1\}\in H$) and $H_2$ consists of the remaining edges of $H$ 
not containing $1$. We relabel the vertices by an increasing injection so that $\bigcup H_1=[k]$ and 
$\bigcup H_2=[l]$. It is clear that 
$H_1$ and $H_2$ are simple and that $k,l\le n-1$. To invert the decomposition, we first select two simple hypergraphs $H_1$ and 
$H_2$ with $\bigcup H_1=[k]$ and $\bigcup H_2=[l]$, which can be done in $h_kh_l$ ways. We relabel their vertices 
and unite the vertex sets so that the set $[2,n]$ arises. 
This can be done in exactly 
$$
{n-1\choose k+l-n+1,\;n-1-k,\;n-1-l}
$$ 
ways by partitioning $[2,n]$ in $k+l-n+1$, $n-1-k$, and $n-1-l$ vertices lying in $C=\bigcup H_1\cap\bigcup H_2$, 
$\bigcup H_2\backslash C$, and 
$\bigcup H_1\backslash C$, respectively. We append to every edge in $H_1$ the new least vertex $1$ and obtain 
a simple hypergraph $H$ with 
$n$ vertices. Finally, the possible addition of $\{1\}$ to $H$ (we always loose the edge $\{1\}$ when decomposing) gives 
two further options, with the exception of $H_1=\emptyset$ when $\{1\}$ must be always added. This explains the factor 2 and 
the subtraction of $h_{n-1}$. The stated recurrence follows. 
\kduk

\noindent
Either of the recurrences 1 and 3 or the explicit formula 2 give
$$
(h_n)_{n\ge 1}=(1,5,109,32297,2147321017,9223372023970362989,\ldots).
$$
This quickly growing sequence is entry A003465 of Sloane \cite{sloa_hand}. 

We proceed to the problem of counting hypergraphs, simple and all, by their weight. The enumeration of all hypergraphs 
$F$ with $i(F)\le 4$ 
in Section 3 shows that $(h_n')_{n\ge 1}=(1,2,7,28,\ldots)$ and $(h_n'')_{n\ge 1}=(1,3,10,41,\ldots)$. We derive 
some formulae and algorithms which produce further terms of these sequences. Recall that a partition 
$\lambda=1^{a_1}2^{a_2}\ldots l^{a_l}$ of $n\in\N$, where $a_i\ge 0$ are integers and $a_l>0$,  
is the decomposition $n=1+1+\cdots+1+2+\cdots+2+\cdots+l+\cdots+l$ with the part $i$ appearing $a_i$ times. 
Thus $\sum_1^l ia_i=n$. We write briefly $\lambda\vdash n$. 
If the hypergraph $H$ has weight $n$ and $a_i$ edges of cardinality $i$, the maximum 
edge cardinality being $l$, then $\lambda=1^{a_1}2^{a_2}\ldots l^{a_l}\vdash n$ and we say that $H$ has {\em edge type\/} 
$\lambda$. We begin with counting hypergraphs with a fixed edge type.

\begin{veta}\label{hranovytyp}
Let $\lambda=1^{a_1}2^{a_2}\ldots l^{a_l}\vdash n$ where $a_l>0$. 
The number of nonisomorphic simple hypergraphs with weight $n$ and edge type 
$\lambda$ is
$$
\sum_{j=l}^n{{j\choose 1}\choose a_1}{ 
{j\choose 2}\choose a_2}\ldots{{j\choose l}\choose a_l}
\sum_{m=j}^n(-1)^{m-j}{m\choose j}
$$
and the number of nonisomorphic hypergraphs with weight $n$ and edge type 
$\lambda$ is
$$
\sum_{j=l}^n{{j\choose 1}+a_1-1\choose a_1}{ 
{j\choose 2}+a_2-1\choose a_2}\ldots{{j\choose l}+a_l-1\choose a_l}
\sum_{m=j}^n(-1)^{m-j}{m\choose j}.
$$
\end{veta}
\duk
Consider the polynomials
$$
W_n=W_n(x_1,x_2,\ldots,x_n)=\sum_{H}\prod_{i=1}^n x_i^{e(i,H)}
$$
where we sum over all simple $H$ with $\bigcup H=[n]$, and $e(i,H)$ is the number of $i$-element edges in $H$. We refine 
the identity from the proof of 1 of Proposition~\ref{hvn} (which corresponds to $x_1=x_2=\cdots=x_n=1$) and obtain
$$
\prod_{i=1}^n (1+x_i)^{{n\choose i}}=\sum_{j=0}^n{n\choose j}W_j.
$$
In terms of exponential generating functions,
$$
\sum_{n\ge 0}{y^n\over n!}\cdot\prod_{i=1}^n (1+x_i)^{{n\choose i}}=
\mathrm{e}^y\cdot\sum_{n\ge 0}{W_ny^n\over n!}.
$$
We invert this relation as in the proof of 2 of Proposition~\ref{hvn} and get
$$
W_n(x_1,\ldots,x_n)=\sum_{j=0}^n(-1)^{n-j}{n\choose j}
\prod_{i=1}^j (1+x_i)^{{j\choose i}}.
$$
The number of nonisomorphic simple hypergraphs $H$ with $i(H)=n$ and edge type $\lambda=1^{a_1}2^{a_2}\ldots l^{a_l}\vdash n$ 
is the coefficient at $x_1^{a_1}\ldots x_l^{a_l}$ in $W_l+W_{l+1}+\cdots+W_n$ which equals 
$$
\sum_{m=l}^n\sum_{j=0}^m(-1)^{m-j}{m\choose 
j}\prod_{i=1}^l{{j\choose i}\choose a_i}=
\sum_{j=l}^n\prod_{i=1}^l{{j\choose i}\choose a_i}\sum_{m=j}^n(-1)^{m-j}{m\choose j}.
$$
The derivation of the second formula is similar, only $W_n$ becomes a power series and $1+x_i$ is replaced by $(1-x_i)^{-1}$ 
because now any $i$-element edge may come in arbitrary many copies. 
\kduk

\noindent
We give for illustration the distribution of hypergraphs with weight $6$ by their edge types. The first 
entry is the number of simple hypergraphs and the second, given only if different, is the number of all hypergraphs:

\bigskip
\begin{center}
\begin{tabular}{l||l|l|l|l|l|l|l}
$\lambda$ & $6^1$ & $1^15^1$ & $2^14^1$ & $1^24^1$ & $3^2$ & $1^12^13^1$
& $1^33^1$ \\
\hline
$\# H$ & $1$ & $11$ & $41$ & $41,50$ & $31,32$ & $239$ & $63,120$ 
\end{tabular}

\bigskip
\begin{tabular}{l|l|l|l}
$2^3$ & $1^22^2$ & $1^42^1$ & $1^6$\\
\hline
$62,75$ & $198,264$ & $41,160$ & $1,32$
\end{tabular}
\end{center}

\bigskip
Collecting the numbers over all edge types, we obtain formulae for the numbers $h_n'$ and $h_n''$.
\begin{dusl}\label{hnishns}
The numbers of nonisomorphic hypergraphs with weight $n$, simple and all, are 
($\lambda=1^{a_1}2^{a_2}\ldots l^{a_l}$ with $a_l>0$)
\begin{eqnarray*}
h_n'&=&\sum_{\lambda\vdash n}\sum_{j=l}^n
\prod_{i=1}^l{{j\choose i}\choose a_i}\sum_{m=j}^n(-1)^{m-j}{m\choose j}\\
h_n''&=&\sum_{\lambda\vdash n}\sum_{j=l}^n
\prod_{i=1}^l{{j\choose i}+a_i-1\choose a_i}\sum_{m=j}^n(-1)^{m-j}{m\choose j}.
\end{eqnarray*}
\end{dusl}

\noindent 
Using these formulae and computer algebra system MAPLE, we have found the following values. 

\bigskip
\begin{center}
\begin{tabular}{l||l|l|l|l|l|l|l|l|l|l}
$n$ &1&2&3&4&5&6&7&8&9&10\\
\hline $h_n'$ & 1 & 2 & 7 & 28 & 134 & 729 & 4408 & 29256 & 210710 & 1633107\\
\hline $h_n''$ & 1 & 3 & 10 & 41 & 192 & 1025 & 6087 & 39754 & 282241 & 2159916
\end{tabular}

\bigskip
\begin{tabular}{l|l}
11 & 12 \\
\hline 13528646 & 119117240\\
\hline 17691161 & 154192692
\end{tabular}
\end{center}

\bigskip
Each of the three formulae in Proposition~\ref{hvn} gives an algorithm that calculates $h_n$ 
in $O(n^c)$ arithmetical operations. In fact, formula 2 requires only $O(n)$ operations.
In contrast, Corollary~\ref{hnishns} gives algorithms that calculate $h_n'$ and $h_n''$ 
in roughly $n^cp(n)$ operations, where $p(n)=|\{\lambda:\ \lambda\vdash n\}|$, which is a superpolynomial number 
because $p(n)\sim(n\cdot 4\sqrt{3})^{-1}\cdot \exp(\pi(2n/3)^{1/2})$ 
as found by Hardy and Ramanujan \cite{hard_rama} (see also Andrews \cite{andr} and Newman \cite{newm62,newm98}). 
From the complexity point of view, Corollary~\ref{hnishns} is much less effective than Proposition~\ref{hvn}. 
On the other hand, it is superior to the trivial way of calculating $h_n'$ and $h_n''$ because these numbers grow 
superexponentially (see Proposition~\ref{asym}) but $p(n)$ is subexponential. The number of operations required by 
Corollary~\ref{hnishns} is therefore still substantially smaller than the number of objects enumerated by $h_n'$ and 
$h_n''$. We show that $h_n'$ and $h_n''$ can be calculated more effectively, again in $O(n^c)$ arithmetical operations, 
by the approach that we used in the recurrence 3 of Proposition~\ref{hvn}. 

To this end we define $h_{n,m,l}'$ to be the number of simple nonisomorphic hypergraphs 
$H$ with $i(H)=n$, $e(H)=m$, and $v(h)=l$. The quantity $h_{n,m,l}''$ is defined similarly for all hypergraphs. 
Obviously, $h_{n,m,l}'=0$ whenever $m>n$ or $l>n$, and the same holds for $h_{n,m,l}''$. Thus, for $n\ge 1$,
$$
h_n'=\sum_{1\le m,l\le n}h_{n,m,l}'\ \mbox{ and }\ h_n''=\sum_{1\le m,l\le n}h_{n,m,l}''.
$$
These sums have $n^2$ summands. To obtain an effective algorithm for calculating $h_n'$ and $h_n''$, it suffices to establish
effective recurrent relations for $h_{n,m,l}'$ and $h_{n,m,l}''$.

\begin{prop}\label{efectrecurr}
Let $T(a,b,c)={a\choose b+c-a,a-b,a-c}$. We have $h_{0,0,0}'=h_{0,0,0}''=1$, $h'_{n,m,l}=h''_{n,m,l}=0$ if $nml=0$ but 
$n+m+l>0$, and, for $n\ge 1$ and $1\le m,l\le n$, 
\begin{eqnarray*}
h_{n,m,l}'&=&\sum_{p=1}^m\sum T(l-1,l_1,l_2)\cdot(h_{n_1,p,l_1}'+h_{n_1,p-1,l_1}') h_{n_2,m-p,l_2}'\\
h_{n,m,l}''&=&\sum_{p=1}^m\sum T(l-1,l_1,l_2)\sum_{q=0}^p h_{n_1,p-q,l_1}''h_{n_2,m-p,l_2}'',
\end{eqnarray*}
where the summation range of the second sum is in both formulae $n_i\ge 0$, $l_i\ge 0$, $n_1+n_2=n-p$, and 
$\max(l_1,l_2)\le l-1\le l_1+l_2$. 
\end{prop}
\duk
We begin with the case of simple hypergraphs.
We decompose any simple hypergraph $H$ with $i(H)=n$, 
$e(H)=m$, and $\bigcup H=[l]$ in the hypergraphs $H_1$ and $H_2$, where $H_1=(E\backslash\{1\}:\ 1\in E\in H)$ and 
$H_1=(E:\ 1\not\in E\in H)$. If $\{1\}\in H$, we remove $\emptyset$ from $H_1$. 
We denote $p=\deg_H(1)$, $i(H_1)=n_1$, $i(H_2)=n_2$, $v(H_1)=l_1$, and $v(H_2)=l_2$. It is clear that $e(H_2)=m-p$
and that the conditions of the second sum are met. If $\{1\}\not\in H$ then $e(H_1)=p$ else $e(H_1)=p-1$. The decomposition
is inverted as in the proof of 3 of Proposition~\ref{hvn}. The cases $\{1\}\not\in H$ and $\{1\}\in H$ are reflected by the 
terms $h_{n_1,p,l_1}'$ and $h_{n_1,p-1,l_1}'$, respectively. The trinomial $T(l-1,l_1,l_2)$ counts the number of ways 
in which the set $[2,l]$ can be written as a union of two sets with $l_1$ and $l_2$ elements. 
We obtain the first 
recurrence. The proof of the recurrence for all hypergraphs is similar, the only difference being that now $\{1\}$ 
may have in $H$ multiplicity $q$, $0\le q\le p=\deg_H(1)$. 
\kduk

\noindent
The recurrences give algoritms that calculate $h_n'$ in $O(n^6)$ operations and $h_n''$ in $O(n^7)$ operations. 

For every rational polynomial $P(m)\in\Q[m]$ it is true that
$$
\sum_{m=0}^{\infty}{P(m)\over m!}=\mathrm{e}\cdot q
$$ 
where $\mathrm{e}=2.71828\ldots$ is Euler number and $q\in\Q$. This follows by expressing $P(m)$ as the 
$\Q$-linear combination in the basis $\{1,m,m(m-1),m(m-1)(m-2),\ldots\}$. One subfamily of this family of identities 
is Dobi\`nski's formula (\cite[Problems 1.9a and 1.13]{lova} and \cite[p. 210]{comt})
$$
\sum_{m=0}^{\infty}{m^n\over m!}=\mathrm{e}\cdot b_n
$$
in which $b_n$ is the $n$-th Bell number (the number of partitions of $[n]$). We present two combinatorial subfamilies 
which are related to hypergraphs.

\begin{dusl}\label{dobinski}
For every $n\in\N$ we have the identities ($\lambda=1^{a_1}2^{a_2}\ldots l^{a_l}$ with $a_l>0$)
\begin{eqnarray*}
\sum_{m=0}^{\infty}{1\over m!}\cdot\sum_{\lambda\vdash n}
\prod_{i=1}^l{{m\choose i}\choose a_i}&=&\mathrm{e}\cdot\sum_{i(H)=n}^*{1\over v(H)!}\\
\sum_{m=0}^{\infty}{1\over m!}\cdot\sum_{\lambda\vdash n}
\prod_{i=1}^l{{m\choose i}+a_i-1\choose a_i}&=&\mathrm{e}\cdot\sum_{i(H)=n}{1\over v(H)!}
\end{eqnarray*}
where $\mathrm{e}=2.71828\ldots$ and the star indicates that the sum is over 
simple hypergraphs $H$ only.
\end{dusl}
\duk
Let $n\in\N$. In the proof of Theorem~\ref{hranovytyp} we used for simple hypergraphs the equation
$$
\sum_{m\ge 0}{y^m\over m!}\cdot\prod_{i=1}^m (1+x_i)^{{m\choose i}}=
\mathrm{e}^y\cdot\sum_{m\ge 0}{W_my^m\over m!}.
$$
The first stated identity now follows by setting $x_i=x^i$, $i\in\N$, comparing the coefficients at 
$x^n$ on both sides, and setting $y=1$. The second identity follows by the same way from the analogous equation for all 
hypergraphs.
\kduk

\noindent
For $n=1,2,3$, and $4$ the factors at $\mathrm{e}$ in the first identity are, respectively, 
$1,1,{11\over 6}$, and ${25\over 8}$, and in the second identity they are $1,2,{23\over 6}$, and ${89\over 8}$. 

It is natural to ask about the asymptotics of $h_n'$ and $h_n''$. We give a simple estimate in terms of the 
Bell numbers $b_n$. 

\begin{prop}\label{asym}
For every $n\in\N$, one has the inequalities
$$
b_n\le h_n'\le h_n''\le 2^{n-1}b_n.
$$
For $n\to\infty$,
$$
\log h_n''=\log b_n+O(n)=n(\log n-\log\log n+O(1))
$$
and the same holds for $h_n'$.
\end{prop} 
\duk
The first two inequalities are trivial. To prove the third inequality, we assign to every hypergraph $H$, where
$i(H)=n$ and $\bigcup H=[m]$ with $m\le n$, a pair $(Q,P)$ of partitions of $[n]$ as follows. 
We set $Q=(I_1,I_2,\dots,I_m)$ where 
$I_1<I_2<\dots<I_m$ are intervals such that $|I_i|=\deg_H(i)$. Thus $Q$ is a partition of $[n]$ into intervals. For every 
$E\in H$ we select a set $A_E\subset [n]$, $|A_E|=|E|$, such that (i) for every $i\in[m]$, 
$A_E\cap I_i\ne\emptyset$ iff $i\in E$ and 
(ii) the sets $A_E$ are mutually disjoint. This can be done and generally in more than one way. We set 
$P=(A_E:\ E\in H)$. It is clear that, regardless of the freedom in selecting $P$, distinct hypergraphs $H$ produce distinct 
pairs $(Q,P)$. The number of pairs $(Q,P)$ does not exceed $2^{n-1}b_n$ because there are exactly $2^{n-1}$ interval 
partitions of $[n]$. Thus we have the inequality $h_n''\le 2^{n-1}b_n$. The logarithmic asymptotics follows from 
the asymptotics of $b_n$ that was found by Moser and Wyman \cite{mose_wyma}, see \cite[Problem 1.9b]{lova} or Odlyzko 
\cite{odly}. 
\kduk

\noindent
It is an interesting question how tight is each of the three above inequlities. The previous argument made no use 
of the fact that the partitions $Q$ and $P$ are ``orthogonal'' in the sense that $|I\cap A|\le 1$ for every $I\in Q$ and 
$A\in P$. Using this, we can narrow the gap in the estimate $b_n\le h_n''\le 2^{n-1}b_n$. We shall treat this topic 
elsewhere.

\end{document}